\documentclass[12pt]{article}

\newcommand{\address}{\noindent
Faculty of Mathematics, Kyushu University 33, \\
Fukuoka 812-8581, Japan. \\
Email: {\tt mkaneko@math.kyushu.u-ac.jp} 
}

\usepackage{amsmath} \usepackage{amssymb} \usepackage{amsfonts}
\usepackage{float}

\def\a{\alpha}  \def\A{{\cal A}}    \def\d{\delta}
\def\disc{{\rm disc}} \def\g{\gamma}  \def\G{\Gamma}
\def\e{\varepsilon}  \def\G{\Gamma}  
  \def\mg{\PSL_2(\Z)}
   
\def \PGL{{\rm PGL}}\def \PSL{{\rm PSL}} 
\def\pf{\noindent{\it Proof. }}  \def\qed{\quad$\square$\medskip } 
  \def\R{\mathbf{R}}
 \def\ta{\tau} \def\uh{\mathfrak{H}}
\def\val{{\rm val}} \def\Z{\mathbf{Z}}

\title{Observations on the ``values'' of the elliptic
modular function $j(\tau)$ at real quadratics
}
\author{Masanobu Kaneko}
\date{}

\begin{document}

\maketitle

\begin{abstract}
We define ``values'' of the elliptic modular $j$-function at
real quadratic irrationalities by using Hecke's hyperbolic
Fourier expansions, and present some observations based 
on numerical experiments. 
\end{abstract}




\section{Introduction}

\begin{quotation}
\noindent ``.... von dem Studium des Verhaltens der 
elliptischen Modulfunktionen in der N\"ahe der nicht-rationalen
Randpunkte noch sehr bemerkenswerte Ergebnisse erwartet werden
k\"onnen, die sowohl f\"ur die Funktionentheorie  wie die Arithmetik
wichtig sein d\"urften.''   (Hecke, Werke S.417)
\end{quotation}
\medskip	

We define the ``value,'' written $\val(w)$,\footnote{Dedekind,  in his
  seminal paper \cite{D} on $j(\tau)$,  used the symbol $\val(\omega)$ for
  $j(\tau)$ (where $\omega$ is a variable in the upper half-plane) and
  called it the ``Valenz.'' We borrow his notation.} of  the elliptic
modular function $j(\ta)$ at each {\it real} quadratic irrationality
$w$ as the constant term of a hyperbolic Fourier
expansion\footnote{Hecke considered this type of expansion for modular
  forms of positive weight \cite{H}.} at $w$. The map $w\mapsto\val(w)$ is
$\PSL_2(\Z)$-invariant and hence assigns to each
$\PSL_2(\Z)$-equivalence class of  real quadratic numbers a certain
(real or complex) number.  We conducted numerical experiments on the
numbers $\val(w)$ and observed the following phenomena, which we find
quite remarkable, though no precise formulation (especially for
(ii) and (iii)) nor proofs have yet been established:
\medskip

\noindent{\bf  Observations}

(i) {\em The minimum  among all real values of $\val(w)$ is  realized
  at $w=(1+\sqrt{5})/2$ (the golden ratio), with
  $\val((1+\sqrt{5})/2)=706.324813540\ldots$. Also, all real values of
  $\val(w)$ lie in the interval $[706.324813540\ldots, 744]$, where
  $744$ is the constant term in the Fourier expansion of $j(\tau)$ at
  the cusp (which is the $\PSL_2(\Z)$-equivalence class of rational
  numbers and $i\infty$).} \smallskip
  
(ii) {\em As the rational approximation of $w$ improves, $\val(w)$
  increases. (See the tables at the  end of
  the paper.)} \smallskip

(iii) {\em The imaginary part of any $\val(w)$ lies in the interval
  $(-1,1)$. Also, the distribution of the imaginary parts of
  $\val(w)$, with the discriminants of $w$ bounded, seems to be peaked
  at 0 and symmetric about this peak.  Furthermore, the  phenomena
  described in (i) and (ii) also hold for the absolute value (or real part) of
  $\val(w)$.  } \bigskip

In this paper, we give a precise definition of $\val(w)$ and then
establish its basic properties, which follow almost
immediately from the definition.  We then describe experiments related to 
Markoff numbers.  This also seems to support the existence of
certain ``Diophantine continuity'' of $\val$ suggested (but not yet well-formulated)
above.   

\section{Definition and basic properties}

Let $w$ be a real quadratic number with discriminant
$\disc(w)=D>0$. Denote by $\Gamma_w$ the stabilizer of $w$ in
$\Gamma=\PSL_2(\Z)$ (with the action being the standard linear
fractional transformation):
$$\Gamma_w:=\left\{\gamma\in\Gamma\, |\, \gamma w=w\right\}.$$
Let
$U_D$ be the group of units of norm 1 in the quadratic order $O_D$ of
the discriminant $D$ and $\e=\e_D^{(1)}$ be a generator of the infinite
cyclic part of $U_D$. Then, if $\gamma=\pm\begin{pmatrix} a & b\\ c & d
\end{pmatrix}\in\Gamma_w$,  we have
$cw^2+(d-a)w-b=0$, and thus the number $cw$ is an 
algebraic integer, and  $$ (a-cw)(a-cw')=a^2-ac(w+w')+c^2ww'=1,$$
that is, $a-cw\in U_D$. Here, $w'$ is the
algebraic conjugate of $w$. It is known that the map
$$
\Gamma_w\ni \gamma=\pm\begin{pmatrix} a & b\\ c & d \end{pmatrix}
\mapsto (a-cw)^2\in U_D^2$$
gives an isomorphism from the group
$\Gamma_w$ to $U_D^2$, which is an infinite cyclic group generated by
$\e^2$. Let $\g_\e$ be the element in $\G_w$ that corresponds to
$\e^2$  under this isomorphism.   For $\gamma\in\Gamma_w$, a
straightforward computation shows $$\frac{\gamma\tau-w}
{\gamma\tau-w'}=(a-cw)^2\cdot\frac{\tau-w}{\tau-w'}$$
and, in particular,
$$\frac{\gamma_\e\tau-w}
{\gamma_\e\tau-w'}=\e^2\frac{\tau-w}{\tau-w'}.$$
Denote by $\d(w)$ the sign of $w-w'$.
Then, if $\ta$ is a variable in the upper half plane $\uh$, we have
$$z:=\d(w)\frac{\ta-w}{\ta-w'}\in \uh,$$
and
$$\ta=\frac{w-\d(w)w'z}{1-\d(w)z}.$$ Let $$j(\ta)=q^{-1}+744+196884q+
21493760q^2+\cdots
\quad (q=e^{2\pi i\ta})$$ be {\it the}
classical elliptic modular function.
It is $\G$-invariant, and hence, by the relations
$$\frac{\g_w\ta-w}{\g_w\ta-w'}=\e^2\d(w)z$$ and $$\g_w\ta=
\frac{w-\e^2\d(w)w'z}{1-\e^2\d(w)z},$$ the function
$$j(\tau)=j\left(\frac{w-\d(w)w'z}{1-\d(w)z}\right)\quad (z\in\uh)$$
is invariant under $z\mapsto\e^2z$. Thus, if we set $z=e^u$, the
function $$j\left(\frac{w-\d(w)w'e^u}{1-\d(w)e^u}\right),$$
which is
holomorphic in the domain $0<\rm{Im}(u)<\pi$, is invariant under the
translation $u\mapsto u+2\log\e$. It therefore has a Fourier
expansion of the form
\begin{equation}\label{Fourier}
j\left(\frac{w-\d(w)w'e^u}{1-\d(w)e^u}\right)=
\sum_{n=-\infty}^\infty a_ne^{2\pi i n
  \frac{u}{2\log\e}}. 
\end{equation}
\medskip

\noindent{\bf Definition} {\em We define the ``value,'' $\val(w)$, 
of $j(\ta)$ at  $w$ as the constant term of the series \eqref{Fourier}:
\begin{equation}\label{a0formula0}
\val(w):=a_0=\frac1{2\log\e}\int_{\sigma_0}^{\sigma_0+2\log\e} j\left(
\frac{w-\d(w)w'e^u}{1-\d(w)e^u}\right)\,du,
\end{equation}
where $\sigma_0$ is any complex number satisfying $0<\rm{Im}(\sigma_0)<\pi$.}

If we set $\sigma_0=\pi i/2-\log\e$ and make the change of variable
$u\mapsto u+\pi i/2$, we have
\begin{equation}\label{a0formula}
\val(w)=\frac1{2\log\e}\int_{-\log\e}^{\log\e} j\left(
\frac{w-\d(w)w'ie^u}{1-\d(w)ie^u}\right)\,du.
\end{equation}
Note that $\val(w)$ is a complex-valued function defined 
{\it only} on the real quadratic irrationalities.\\

\noindent{\bf Proposition} {\em The ``value'' function $\val(w)$ possesses the 
following properties.
 
  1) If $w$ and $w_1$ are $\G$-equivalent, then
  $\val(w)=\val(w_1)$. 

2) $\val(w)=\val(w').$

3) $\overline{\val(w)}=\val(-w').$ }

\pf 1)  Let $w_1=(aw+b)/(cw+d)$, with $\pm\begin{pmatrix} a&b\\c&d 
\end{pmatrix}\in\G$. Because $j(\ta)$ is $\G$-invariant, we have 
\begin{eqnarray*}
j\left(\frac{w-\d(w)w'e^u}{1-\d(w)e^u}\right)&=&
j\left(\frac{a\frac{w-\d(w)w'e^u}{1-\d(w)e^u}+b}
{c\frac{w-\d(w)w'e^u}{1-\d(w)e^u}+d}\right)
=j\left(\frac{aw+b-\d(w)(aw'+b)e^u}
{cw+d-\d(w)(cw'+d)e^u}\right)\\
&=&j\left(\frac{\frac{aw+b}{cw+d}-\d(w)\frac{aw'+b}{cw+d}e^u}
{1-\d(w)\frac{cw'+d}{cw+d}e^u}\right)
=j\left(\frac{w_1-\d(w)w_1'\eta e^u}{1-\d(w)\eta e^u}\right)\\
&=&j\left(\frac{w_1-\d(w_1)
w_1'{\rm sgn}(\eta)\eta e^u}{1-\d(w_1){\rm sgn}(\eta)\eta e^u}\right),
\end{eqnarray*}
where $\eta=(cw'+d)/(cw+d)$ and we have used $\d(w)=\d(w_1){\rm sgn}(\eta)$
[because $w_1-w_1'=(w-w')/((cw+d)(cw'+d))=(w-w')\eta/(cw'+d)^2$].
Therefore, from \eqref{a0formula0}, we obtain 
\[ \val(w)=\frac1{2\log\e}\int_{\sigma_0}^{\sigma_0+2\log\e}
j\left(\frac{w_1-\d(w_1) w_1{\rm sgn}(\eta)\eta e^u}{1-\d(w_1){\rm
      sgn}(\eta)\eta e^u}\right)\,du .\]
Then, because  ${\rm
  sgn}(\eta)\eta>0$, we can make the change of variable $u\mapsto
u-\log({\rm sgn}(\eta)\eta)$, and we conclude that $\val(w)=\val(w_1).$

2)  Changing $u$ to $-u$ in \eqref{a0formula} and using the relation
$\d(w')=-\d(w)$, we have

 \begin{eqnarray*}
 \val(w)&=& \frac1{2\log\e}\int_{-\log\e}^{\log\e} j\left(
 \frac{w-\d(w)w'ie^{-u}}{1-\d(w)ie^{-u}}\right)\,du 
 = \frac1{2\log\e}\int_{-\log\e}^{\log\e} j\left(
 \frac{w'+\d(w)wie^u}{1+\d(w)ie^u}\right)\,du \\
 &=& \frac1{2\log\e}\int_{-\log\e}^{\log\e} j\left(
 \frac{w'-\d(w')(w')'ie^u}{1-\d(w')ie^u}\right)\,du = \val(w'). 
 \end{eqnarray*}

3) By \eqref{a0formula}, we have \begin{eqnarray*}
\overline{\val(w)}&=& \frac1{2\log\e}\int_{-\log\e}^{\log\e}
\overline{j\left(\frac{w-\d(w)w'ie^u}{1-\d(w)ie^u}\right)}\,du
=\frac1{2\log\e}\int_{-\log\e}^{\log\e}j\left(-\frac
{w+\d(w)w'ie^u}{1+\d(w)ie^u}\right)\,du\\
&=&\frac1{2\log\e}\int_{-\log\e}^{\log\e}j\left(\frac
{-w'+\d(w)wie^{-u}}{1-\d(w)ie^{-u}}\right)\,du\\
&=&\frac1{2\log\e}\int_{-\log\e}^{\log\e}j\left(\frac
{-w'-\d(-w')(-w')'ie^{-u}}{1-\d(-w')ie^{-u}}\right)\,du
=\val(-w').
\end{eqnarray*}
\qed

\noindent{\it Remark.} The  invariance in 1) does not hold in general 
for other coefficients $a_n=a_n(w)\  (n\ne0)$. The general transformation 
formula is similarly deduced and reads
\[  a_n\left(\frac{aw+b}{cw+d}\right)=\left\vert\frac{cw'+d}{cw+d}\right\vert^{-\pi i n/\log\e}
a_n(w). \]

\noindent{\bf Corollary} {\em
  1) Suppose $\disc(w)=D$ and let $\e_D$ be the fundamental unit of
the order $O_D$. Then, if $N(\e_D):=\e_D\e_D'=-1$, we always have $\val(w)\in\R$.

2) If $w$ and $-w'$ are $\G$-equivalent, then $\val(w)\in\R$.}\\

\pf 1) In this case, $w$ and $-w$ are $\G$-equivalent, and thus, by
applying 3), 2) and 1) of the Proposition in turn, we obtain
$$\overline{\val(w)}=\val(-w')=\val(-w)=\val(w).$$

2) This follows from 3) and 1) of the Proposition. \qed

We denote by $\A$ the class in the narrow ideal class group $Cl^+(D)$
to which the ideal corresponding to $w$ belongs. By Proposition,
$\val(w)$ depends  only on the class $\A$. (With this in mind,
we may write $\val(\A)$.) The class corresponding to
$-w'$ is ${\A}^{-1}$, and hence the $\G$-equivalence of $w$ and $-w'$ implies
$\A^2=1$ and vice versa.  Hence, the assertion 2) in the corollary
says that the value $\val(\A)$ is real if $\A^2=1$.\\

\noindent{\it Remark.} Numerical computations reveal that 
not all $\val(w)$ are real. \\

We give three examples. \\

\noindent{\bf Example 1.} The minimal discriminant for which there
appears a non-real value is $D=136$. The wide class number $h$ is 2, and
the  narrow one $h^+$ is 4.  A representative of the $\G$-equivalence class
of numbers of discriminant 136 is given by 
\[ \sqrt{34},\
\frac{-4+\sqrt{34}}{18},\ \frac{-1+\sqrt{34}}{11}, \
\frac{1+\sqrt{34}}{11}, \]
and these are grouped into two  wide (${\rm
  PGL}_2(\Z)$-equivalence) classes:
\[ \{\sqrt{34},\,\frac{-4+\sqrt{34}}{18}\},\
\{\frac{-1+\sqrt{34}}{11},\,\frac{1+\sqrt{34}}{11}\}. \]
The narrow class group $Cl^+(136)$ is 
isomorphic to $\Z/4\Z$ and is generated by the class corresponding 
to $(-1+\sqrt{34})/11$.  The values of $\val$ at this generator and
its inverse $(1+\sqrt{34})/11$ (this is also a generator of $Cl^+(136)$) are
computed as
\begin{eqnarray*}
\val\left(\frac{-1+\sqrt{34}}{11}\right)&=&710.600451944002489\ldots - 
0.5197938281961062\ldots i,\\
\val\left(\frac{1+\sqrt{34}}{11}\right)&=&710.600451944002489\ldots +
0.5197938281961062\ldots i,
\end{eqnarray*}
the two being conjugate with each other as follows from Proposition 3).

The values at other two points are
\[ \val(\sqrt{34})=\val\left(\frac{-4+\sqrt{34}}{18}\right)=
720.29003500445066239,\ldots\]
values being identical because $(-4+\sqrt{34})/18$ and $-\sqrt{34}=(\sqrt{34})'$
are $\mg$-equivalent.\\

\noindent{\bf Example 2.} Consider the discriminant
$D=145$. In this case, we have $h=h^+=4$. 
As representative numbers, we may choose 
\[ \frac{1+\sqrt{145}}2,\
\frac{1+\sqrt{145}}6,\ \frac{-5+\sqrt{145}}{12},\
\frac{7+\sqrt{145}}{16}. \]
By Corollary 1) we know all values of $\val$ at these
points are real. Numerically, they are given as
\begin{eqnarray*}
\val\left(\frac{1+\sqrt{145}}2\right)&=&
720.484777347009813\ldots,\\
\val\left(\frac{1+\sqrt{145}}6\right)&=&
715.729503630174741\ldots,\\
\val\left(\frac{-5+\sqrt{145}}{12}\right)&=&
708.568357453922648\ldots,\\
\val\left(\frac{7+\sqrt{145}}{16}\right)&=&
715.729503630174741\ldots.
\end{eqnarray*}
The class group is isomorphic to $\Z/4\Z$.  This is seen from
the fact that, for $w_1= (1+\sqrt{145})/6$, $-w_1'$ is not equivalent 
to $w_1$ but equivalent to $w_2=(7+\sqrt{145})/16$.  Hence
$\val(w_1)=\val(w_2)$. \\

\noindent{\bf Example 3.} Consider $D=520$. In this case, again, 
we have $h=h^+=4$. 
As representative numbers, we may choose 
\[ \sqrt{130},\ 
\frac{-1+\sqrt{130}}3,\ \frac{-3+\sqrt{130}}{11},\ 
\frac{-5+\sqrt{130}}{15},\]  and whose ``values'' are
given numerically by
\begin{eqnarray*}
\val(\sqrt{130})&=&
721.700344576590835\ldots,\\
\val\left(\frac{-1+\sqrt{130}}3\right)&=&
719.032996230455907\ldots,\\
\val\left(\frac{-3+\sqrt{130}}{11}\right)&=&
713.022954982182920\ldots,\\
\val\left(\frac{-5+\sqrt{130}}{15}\right)&=&
716.888481219718920\ldots.
\end{eqnarray*}
In this case, the class group is isomorphic to 
$\Z/2\Z\times\Z/2\Z$, and all values appear to be distinct. \\

\section{Experiments related to Markoff numbers} 

First let us recall Markoff's theory. The classical theorem of 
Hurwitz asserts that, for any real irrational 
number $\alpha$, there exist infinitely many rational 
numbers $p/q$ that satisfy
$$
\left| \a-\frac{p}{q}\right|<\frac1{\sqrt{5} q^2}.$$
The constant
$1/\sqrt{5}$ is best possible. But if we exclude as $\a$ the
numbers which are
$\PGL_2(\Z)$-equivalent to the golden ratio $(1+\sqrt{5})/2$, the
constant $1/\sqrt{5}$ improves to  $1/\sqrt{8}$. If we also exclude
the numbers  which are $\PGL_2(\Z)$-equivalent to  $\sqrt{2}$, then we
can take $5/\sqrt{221}$ as the constant. In general, this continues as
follows. There is an infinite sequence of integers called Markoff
numbers, $$\{m_i\}_{i=1}^\infty=\{1,2,5,13,29,34,89,169,194,233,\ldots\},$$ and
associated quadratic irrationalities $\theta_i$ and monotonically
increasing $L_i$ whose limit is 3,  with the following property: ``For
any $i$, if the number $\a$ is not $\PGL_2(\Z)$-equivalent to any of
$\theta_1, \theta_2,\ldots,\theta_{i-1}$, then there exist infinitely
many rational numbers $p/q$ that satisfy
$$
\left| \a-\frac{p}{q}\right|<\frac1{L_i q^2}.$$
Explicitly, the Markoff numbers $m_i$ appear as solutions of the
diophantine equation 
\begin{equation}\label{mar-eq}
x^2+y^2+z^2=3xyz, 
\end{equation}
and
\begin{equation}\label{theta} 
L_i=\sqrt{9-4/m_i^2},\quad \theta_i=\frac{
  -3m_i+2k_i+\sqrt{9m_i^2-4}}{2m_i},\end{equation}
where $k_i$ is an integer that
satisfies $a_ik_i\equiv b_i\pmod{m_i}$ and here $(a_i,b_i,m_i)$ is a
solution of equation \eqref{mar-eq} with $m_i$ maximal.
If $(p,q,r)$ is a solution of \eqref{mar-eq}, then $(p,q,3pq-r)$ and $(p,r,3pr-q)$
are also solutions. This gives to all solutions a structure of 
tree, and we can arrange Markoff numbers like the picture below. 

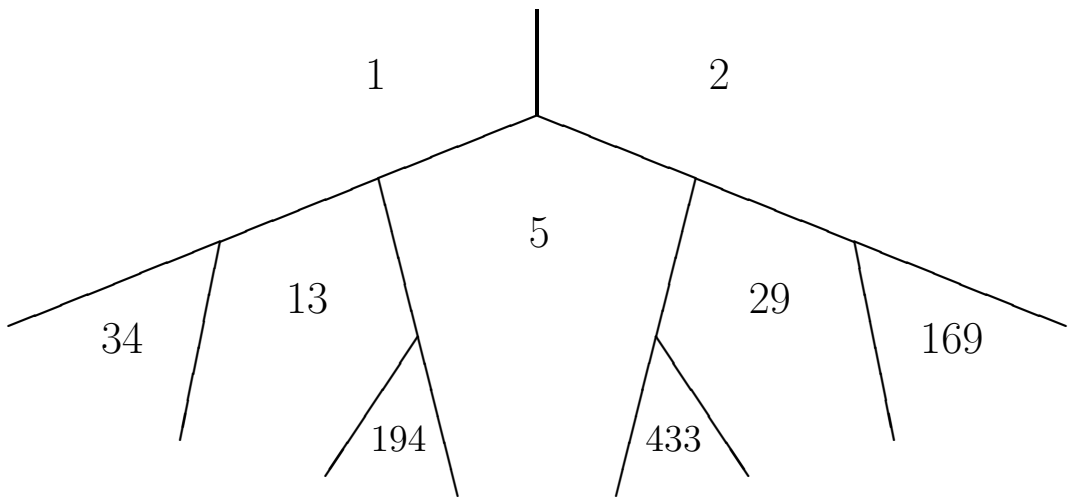
\begin{figure}[H]\caption{The tree of Markoff numbers}
\begin{picture}(200,200)
\thicklines
\put(220,150){\line(0,1){40}}
\put(220,150){\line(5,-2){200}}
\put(220,150){\line(-5,-2){200}}
\put(160,126){\line(1,-4){30}}
\put(280,126){\line(-1,-4){30}}
\put(100,102){\line(-1,-5){15}}
\put(340,102){\line(1,-5){15}}
\put(175,66){\line(-2,-3){35}}
\put(265,66){\line(2,-3){35}}
\put(155,160){\Large{$1$}}
\put(285,160){\Large{$2$}}
\put(217,100){\Large{$5$}}
\put(125,75){\Large{$13$}}
\put(300,75){\Large{$29$}}
\put(55,60){\Large{$34$}}
\put(365,60){\Large{$169$}}
\put(157,23){\large{$194$}}
\put(261,23){\large{$433$}}
\end{picture}
\end{figure}

We computed several values of $\val(\theta_i)$, and observed
the following. \medskip	

\noindent{\bf Observation}  (iv)  {\em Only real values are $$\val(\theta_1)=\val\bigl(\frac{-1+\sqrt{5}}2\bigr)
=706.32481354\ldots $$and $$\val(\theta_2)=\val(-1+\sqrt{2})=709.89289091
\ldots.$$ No other values $\val(\theta_i)\ (i\ge3)$ seem to be real. }\\ 

Note that in Markoff's theory only $\PGL_2(\Z)$-equivalence
class is relevant, but we need $\mg$-equivalence to distinguish
non-real $\val(\theta_i)$ and its conjugate. Here, the order of $(a_i,b_i)$ in
the definition of $\theta_i$ in \eqref{theta} becomes relevant. We introduce
the following refinement.

Let $(a,b,m_i)$ be the Markoff triple associated to the $i$th Markoff number
$m_i$ and assume the order of $a$ and $b$ is so chosen 
that their positions in the tree is like 

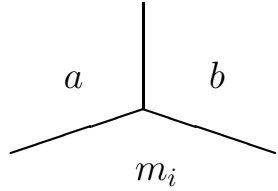
\begin{figure}[H]
\begin{picture}(100,80)
\thicklines
\put(170,45){\large{$a$}}
\put(200,37){\line(0,1){40}}
\put(225,45){\large{$b$}}
\put(200,37){\line(3,-1){50}}
\put(197,10){\large{$m_i$}}
\put(200,37){\line(-3,-1){50}}
\end{picture}
\end{figure}

\noindent ($b$ is on the right of $a$, so, $(13,5,194)$ for $194$, $(5,29,433)$
for $433$ etc.)
Define two numbers $\theta_{i,1}$ and $\theta_{i,2}$ by 
$$\theta_{i,1}=\frac{-3m_i+2k_{i,1}+\sqrt{9m_i^2-4}}{2m_i}\quad\text{and}\quad
\theta_{i,2}=\frac{  -3m_i+2k_{i,2}+\sqrt{9m_i^2-4}}{2m_i}$$
with $k_{i,1}$ and $k_{i,2}$ being integers that satisfy
$$ak_{i,1}\equiv b\pmod{m_i}\quad\text{and}\quad
 bk_{i,2}\equiv a\pmod{m_i}$$ respectively. \\
 
\noindent{\bf Observation} (v) {\em The imaginary part of $\val(\theta_{i,1})$ 
(resp. $\val(\theta_{i,2})$) is always positive (resp. negative).} \\

\noindent{\bf Observation} (vi) {\em Suppose three Markoff numbers $m,\,m',\,m''$ are 
in the position like 

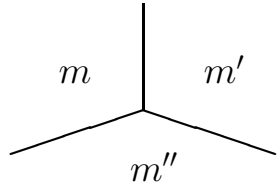
\begin{figure}[H]
\begin{picture}(100,80)
\thicklines
\put(168,47){\large{$m$}}
\put(200,37){\line(0,1){40}}
\put(223,47){\large{$m'$}}
\put(200,37){\line(3,-1){50}}
\put(195,10){\large{$m''$}}
\put(200,37){\line(-3,-1){50}}
\end{picture}
\end{figure}

\noindent in the Markoff tree, and let $\theta_1,\,\theta_2,\,\theta_1',\,\theta_2',\,\theta_1'',\,\theta_2''$ be 
the associated (refined) quadratic numbers. Then, for $j=1,2$, both the real and the imaginary parts 
of $\theta_j''$ lie between those of $\theta_j$ and $\theta_j'$ (the case of $m=1,\,m'=2,\,m''=5$ is exceptional, where the imaginary
parts of $\val(\theta_j)$ and $\val(\theta_j')$ are both 0, while the real part of $\val(\theta_j'')$
is indeed in between those of $\val(\theta_j)$ and $\val(\theta_j')$).}

Hence, all real parts of $\val(\theta_{i,j})\ (j=1,2)$, conjecturally, lie
in the interval $$[706.3248135\ldots, 709.8928909\ldots]$$
and imaginary parts in $$[-0.2670397\ldots, 0.2670397\ldots],$$
where $0.2670397\ldots$ is the imaginary part of 
$$\val(\theta_{3,1})=\val((-11+\sqrt{221})/10)=708.90991972\ldots
+0.267039735\ldots i.$$
\medskip

Choose any Markoff number $m$. This determines a connected unbounded
region $R$ in the tree. 
If we trace the edges of $R$ downward, we obtain the sequence of Markoff numbers
associated to the neighboring region with respect to 
those edges. Let 
$$n_1^L,n_2^L, n_3^L,\ldots\quad\text{and}\quad n_1^R,n_2^R, n_3^R,\ldots$$
be those sequences corresponding to the left and the right edges respectively.
(When $m=1$ (resp. $m=2$), only the sequence $\{n_k^R\}$ (resp. $\{n_k^L\}$) occur.) \\

\noindent{\bf Observation} (vii)  {\em Let $\theta_1^{(m)}$ and $\theta_2^{(m)}$ be the 
Markoff irrationalities associated to $m$ as explained
above (by fixing the order of $a$ and $b$ in the triple $(a,b,m)$), and similarly
$\theta_{k,j}^L\ (j=1,2)$ (resp. $\theta_{k,j}^R\ (j=1,2)$ ) the irrationalities 
associated to $n_k^L$ (resp. $n_k^R$). Then, we surmise
$$\lim_{k\to\infty}\val(\theta_{k,1}^R)=\val(\theta_1^{(m)})\quad\text{and}\quad
\lim_{k\to\infty}\val(\theta_{k,2}^L)=\val(\theta_2^{(m)}).$$ }
\medskip	

Below, we repeat the observations made at the beginning of the paper,
in the form of several questions:
\begin{itemize}
\item Is $\val((1+\sqrt{5})/2)=706.32481354081\ldots$ minimal
  (in absolute value) among all the values of $j(\ta)$ at real
  quadratics? Do all real values of $\val(w)$, or all absolute values 
  or real parts of $\val(w)$, lie in the interval
  $[706.32481354081\ldots, 744]$? If this is the case, is $744$
  the best possible upper bound?
\item Does $\val(w)$ possess some information concerning the Diophantine
  approximation of $w$?  For instance, does $\val(w)$ increase as the
  rational approximation of $w$ improves?  
\item Does the imaginary part of $\val(w)$ always lie in the 
interval $(-1,1)$? What is the distribution of the imaginary 
parts?
\end{itemize}

\noindent{\bf Problem} {\it Formulate rigorous statements and
  find proofs of them that answer all of these questions and,  above
  all, find an arithmetic meaning of $\val(w)$.}
\medskip

\noindent{\it Remark.}  1)  Concerning the nature of the value  $\val(w)$, 
numerical experiments suggest that it is very unlikely that $\val(w)$ is
itself an algebraic number.  
The author has spent a fair amount of time, using ``lindep'' or 
``algdep'' facilities of Pari-GP, or ``Plouffe's inverter'' website,
to see if any multiplicative combination of $\val(w)$, $\log\varepsilon$, $\pi$ etc. becomes algebraic, but all in vain so far.

2)  Recent work of Duke, Imamo$\bar{{\rm g}}$lu and T\'oth \cite{DIT} reveals 
that the ``trace'' of $\val(w)$ appears as the Fourier coefficient of a weakly 
harmonic modular forms of weight $1/2$. It would be an important
problem to understand our observations in light of their results.\\

At the end of the paper, we present some tables of values of $\val(w)$.  The computations
were carried out using Mathematica.  

We denote by $[b_1,b_2,\ldots,b_n]$ a purely periodic (ordinary) continued fraction of  period 
length $n$. 
For example, we have $[1]=(1+\sqrt{5})/2,\ [2,1]=1+\sqrt{3}$, etc.  
The fundamental unit of norm 1 (a generator of $U_D$ in \S2) 
of the order $O_D$ of discriminant $D$ is denoted by $\e$. \\

\renewcommand{\baselinestretch}{1.2}
\renewcommand{\arraystretch}{1.1}

\begin{table}[H] \caption{Values of $\val(w)$ at $w=[n]$.}\nonumber
  \begin{center}
  \begin{tabular}{c|c|c|c}
  $w$ & $D$ & $\val(w)$ & $\log\e$ \\
\hline
$[1]$ & $5$ & $706.3248135408125820559603\ldots$ & $0.9624236501192\ldots$ \\
$[2]$ & $8$ & $709.8928909199123368059253\ldots$ & $1.7627471740390\ldots$ \\
$[3]$ & $13$ &$713.2227192129106375260272\ldots$ & $2.3895264345742\ldots$ \\
$[4]$ & $20$ &$715.8658310509644567882877\ldots$ & $2.8872709503576\ldots$ \\
$[5]$ & $29$ &$717.9165510885627097946754\ldots$ & $3.2944622927421\ldots$ \\
$[6]$ & $40$ &$719.5292195149241565812037\ldots$ & $3.6368929184641\ldots$ \\
$[7]$ & $53$ &$720.8247553829016929089184\ldots$ & $3.9314409432993\ldots$ \\
$[8]$ & $68$ &$721.8878326202869588905005\ldots$ & $4.1894250945222\ldots$ \\
$[9]$ & $85$ &$722.7768914565219262830724\ldots$ & $4.4186954172306\ldots$ \\
$[10]$&$104$ &$723.5327700907464960378584\ldots$ & $4.6248766825455\ldots$ \\
$[20]$&$404$ &$727.6296000047325464824629\ldots$ & $5.9964459005959\ldots$ \\
$[30]$&$904$ &$729.4314438625732480951697\ldots$ & $6.8046132909611\ldots$ \\
$[50]$&$2504$ &$731.2426027524741005593885\ldots$ & $7.8248455312825\ldots$ \\
$[100]$&$10004$ &$733.1113065597372736130899\ldots$ & $9.2105403419828\ldots$ 
\end{tabular}
\end{center}
\end{table}

\begin{table}[H] \caption{Values of $\val(w)$ at $w=[n,1]$.}
  \begin{center}
  \begin{tabular}{c|c|c|c}
  $w$ & $D$ & $\val(w)$ & $\log\e$ \\
\hline
$[2,1]$ & $12$ & $709.7923590080320102702826\ldots$ & $1.3169578969248\ldots$ \\
$[3,1]$ & $21$ & $713.2461372719263413372589\ldots$ & $1.5667992369724\ldots$ \\
$[4,1]$ & $32$ & $715.8764861800141880351424\ldots$ & $1.7627471740390\ldots$ \\
$[5,1]$ & $45$ & $717.8834096374473486546884\ldots$ & $1.9248473002384\ldots$ \\
$[6,1]$ & $60$ & $719.4559616552358003854302\ldots$ & $2.0634370688955\ldots$ \\
$[7,1]$ & $77$ & $720.7215682962489544550810\ldots$ & $2.1846437916051\ldots$ \\
$[8,1]$ & $96$ & $721.7640368038035489169855\ldots$ & $2.2924316695611\ldots$ \\
$[9,1]$ & $117$ &$722.6396242176524465181309\ldots$ & $2.3895264345742\ldots$ \\
$[10,1]$ & $140$&$723.3871879544329222875427\ldots$ & $2.4778887302884\ldots$ \\
$[20,1]$ & $480$&$727.4935574326730521838984\ldots$ & $3.0889699048446\ldots$ \\
$[30,1]$ & $1020$&$729.3240631373043636667693\ldots$ & $3.4647579066758\ldots$ \\
$[50,1]$ & $2700$&$731.1703417153756088105933\ldots$ & $3.9508736907744\ldots$ \\
$[100,1]$ & $10400$&$733.0728964687665155522285\ldots$ & $4.6248766825455\ldots$ 
\end{tabular}
\end{center}
\end{table}

\begin{table}[H] \caption{Values of $\val(w)$ at $w=[n,2]$.}
  \begin{center}
  \begin{tabular}{c|c|c|c}
  $w$ & $D$ & $\val(w)$ & $\log\e$ \\
\hline
$[3,2]$ & $60$ & $711.9275163995819056553017\ldots$ & $2.0634370688955\ldots$ \\
$[4,2]$ & $24$ & $713.8258642873420364918902\ldots$ & $2.2924316695611\ldots$ \\
$[5,2]$ & $140$ & $715.4007874465895012696492\ldots$ & $2.4778887302884\ldots$ \\
$[6,2]$ & $48$ & $716.6952844238825705424260\ldots$ & $2.6339157938496\ldots$ \\
$[7,2]$ & $252$ & $717.7711201642989402376217\ldots$ & $2.7686593833135\ldots$ \\
$[8,2]$ & $80$ & $718.6786015779022038417819\ldots$ & $2.8872709503576\ldots$ \\
$[9,2]$ & $396$ &$719.4552346952050033894397\ldots$ & $2.9932228461263\ldots$ \\
$[10,2]$ & $120$&$720.1286213941960093536607\ldots$ & $3.0889699048446\ldots$ 
\end{tabular}
\end{center}
\end{table}

\begin{table}[H] \caption{Values of $\val(w)$ at $w=[2,1,\ldots,1]$.}
  \begin{center}
  \begin{tabular}{c|c|c|c}
  $w$ & $D$ & $\val(w)$ & $\log\e$ \\
\hline
$[2]$ & $8$ &    $709.8928909199123368059253\ldots$ & $1.7627471740390\ldots$ \\
$[2,1]$ & $12$ & $709.7923590080320102702826\ldots$ & $1.3169578969248\ldots$ \\
$[2,1,1]$ &$40$& $708.5134481348921906198907\ldots$ & $3.6368929184641\ldots$ \\
$[2,1,1,1]$ & $96$ & $708.1560508416661547689422\ldots$ & $2.2924316695611\ldots$ \\
$[2,1,1,1,1]$ & $260$ & $707.8064656210238322953785\ldots$ & $5.5529445614474\ldots$ \\
$[2,1,1,1,1,1]$ & $672$ & $707.5978542380262638805993\ldots$ & $3.2566139548000\ldots$ \\
$[2,1,1,1,1,1,1]$ & $1768$& $707.4305612244349322611838\ldots$ & $7.4764720605230\ldots$ 
\end{tabular}
\end{center}
\end{table}

\begin{table}[H] \caption{Values of $\val(w)$ at $w=[3,1,\ldots,1]$.}
 \begin{center}
  \begin{tabular}{c|c|c|c}
  $w$ & $D$ & $\val(w)$ & $\log\e$ \\
\hline
$[3]$ & $13$ & $713.2227192129106375260272\ldots$ & $2.3895264345742\ldots$ \\
$[3,1]$ & $21$ & $713.2461372719263413372589\ldots$ & $1.5667992369724\ldots$ \\
$[3,1,1]$ &$17$& $711.0460844096550318879502\ldots$ & $4.1894250945222\ldots$ \\
$[3,1,1,1]$ & $165$& $710.3366093961225252087583\ldots$ & $2.5589789770286\ldots$ \\
$[3,1,1,1,1]$ & $445$ & $709.6475354979192849968007\ldots$ & $6.0935646744388\ldots$ \\
$[3,1,1,1,1,1]$ & $288$ & $709.2118541585357042188756\ldots$ & $3.5254943480781\ldots$ \\
$[3,1,1,1,1,1,1]$ & $3029$& $708.8593091672155721790085\ldots$ & $8.0153271998839\ldots$ 
\end{tabular}
\end{center}
\end{table}

\begin{table}[H]  \caption{First several non-real values.}
 \begin{center}
  \begin{tabular}{c|c|c}
  $w$ & $D$ & $\val(w)$  \\
\hline
$(12+\sqrt{34})/11$ & $136$ & $710.60045194400248945\ldots+0.51979382819610620\ldots i$ \\
 $(10+\sqrt{34})/11$ & $136$ & $710.60045194400248945\ldots-0.51979382819610620\ldots i$ \\
$(33+\sqrt{205})/34$ & $205$ & $714.16034018225715592\ldots+0.75363913959038068\ldots i$ \\
$(25+\sqrt{205})/30$ & $205$ & $714.16034018225715592\ldots-0.75363913959038068\ldots i$ \\
$(21+\sqrt{221})/22$ & $221$ & $708.90991972070874730\ldots+0.26703973546028996\ldots i$ \\
$(23+\sqrt{221})/22$ & $221$ & $708.90991972070874730\ldots-0.26703973546028996\ldots i$ \\
$(47+\sqrt{305})/56$ & $305$ & $716.13898693848579303\ldots+0.82184193359696810\ldots i$ \\
$(35+\sqrt{305})/46$ & $305$ & $716.13898693848579303\ldots-0.82184193359696810\ldots i$ \\
$(23+\sqrt{79})/25$ & $316$ & $712.65948582687702503\ldots+0.32545553768732463\ldots i$ \\
 $(13+\sqrt{79})/15$ & $316$ & $712.65948582687702503\ldots-0.32545553768732463\ldots i$ \\
  $(17+\sqrt{79})/15$ & $316$ & $712.65948582687702503\ldots+0.32545553768732463\ldots i$ \\
 $(17+\sqrt{79})/21$ & $316$ & $712.65948582687702503\ldots-0.32545553768732463\ldots i$ \\

\end{tabular}
\end{center}
\end{table}

\begin{table}[H]  \caption{First several values at Markoff irrationalities.}
 \begin{center}
  \begin{tabular}{c|c|c|c}
  $i$ & $m_i$ & $\theta_{i,1}$ & $\val(\theta_{i,1})$  \\
\hline
$1$ & $1$ & $(-3 + \sqrt {5})/2  $ & $706.32481354081258205\ldots$ \\
 $2$ & $2$ & $ -1 + \sqrt {2} $ &     $709.89289091991233680\ldots$ \\
$3$ & $5$ & $ (-11 + \sqrt {221})/10 $ &           $708.909919720708747\ldots
+0.267039735460289\ldots i$ \\
$4$ & $13$ & $ (-29 + \sqrt {1517})/26 $ &        $708.257588242846779\ldots
+0.228635826664936\ldots i$ \\
$5$ & $29$ & $(-63 + \sqrt {7565})/58  $ &        $709.302611667387656\ldots
+0.165196473942199\ldots i$ \\
$6$ & $34$ & $ (-19 + 5\sqrt {26})/17 $ &          $707.858372382696744\ldots
+0.184765335383899\ldots i$ \\
$7$ & $89$ & $ (-199 + \sqrt {71285})/178 $ &   $707.594565998876317\ldots
+0.153386774906169\ldots i$ \\
$8$ & $169$ & $(-367 + \sqrt {257045})/338 $ & $709.469768024657232\ldots+0.118518079083046\ldots i$ \\
$9$ & $194$ & $ (-108 + \sqrt {21170})/97 $ &   $708.534665666479421\ldots
+0.245013213468323\ldots i$ \\
 $10$ & $233$ & $(-521 + \sqrt {488597})/466$&$707.408028846873175\ldots
 +0.130903420887032\ldots i$ \\
\end{tabular}
\end{center}
\end{table}

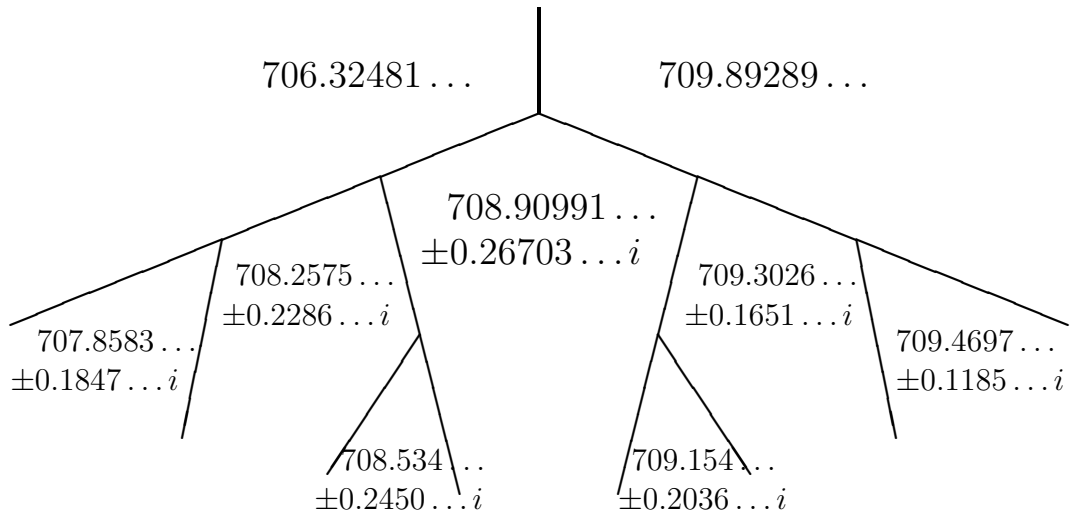
\begin{figure}[H]\caption{Values in the Markoff tree}
\begin{picture}(200,200)
\thicklines
\put(220,150){\line(0,1){40}}
\put(220,150){\line(5,-2){200}}
\put(220,150){\line(-5,-2){200}}
\put(160,126){\line(1,-4){30}}
\put(280,126){\line(-1,-4){30}}
\put(100,102){\line(-1,-5){15}}
\put(340,102){\line(1,-5){15}}
\put(175,66){\line(-2,-3){35}}
\put(265,66){\line(2,-3){35}}
\put(115,160){\large{$706.32481\ldots$}}
\put(265,160){\large{$709.89289\ldots$}}
\put(185,110){\large{$708.90991\ldots$}}
\put(175,93){\large{$\pm 0.26703\ldots i$}}
\put(105,85){{$708.2575\ldots$}}
\put(100,70){{$\pm 0.2286\ldots i$}}
\put(280,85){$709.3026\ldots$}
\put(275,70){$\pm 0.1651\ldots i$}
\put(30,60){$707.8583\ldots$}
\put(20,45){$\pm 0.1847\ldots i$}
\put(355,60){$709.4697\ldots$}
\put(355,45){$\pm 0.1185\ldots i$}
\put(145,15){$708.534\ldots$}
\put(135,0){$\pm 0.2450\ldots i$}
\put(255,15){$709.154\ldots$}
\put(250,0){$\pm 0.2036\ldots i$}
\end{picture}
\end{figure}

\address

\end{document}